\documentclass[12pt]{article}
\usepackage{amsmath, amssymb, amsthm, mathrsfs}
\usepackage[hidelinks]{hyperref}
\usepackage[margin=1in]{geometry}

\theoremstyle{plain}
\newtheorem{theorem}{Theorem}[section]
\newtheorem{lemma}[theorem]{Lemma}
\newtheorem{proposition}[theorem]{Proposition}

\newtheorem{fact}[theorem]{Fact}

\theoremstyle{definition}
\newtheorem{definition}[theorem]{Definition}

\theoremstyle{remark}
\newtheorem{remark}[theorem]{Remark}

\newcommand{\cf}{\mathrm{cf}}
\newcommand{\supp}{\mathrm{supp}}
\newcommand{\ZF}{\mathrm{ZF}}

\newcommand{\DC}{\mathrm{DC}}

\newcommand{\HS}{\mathrm{HS}}

\DeclareMathOperator{\dom}{dom}

\DeclareMathOperator{\sym}{sym}

\newenvironment{acknowledgments}{\par\medskip\noindent\textbf{Acknowledgments.}\ }{\par\medskip}

\title{Limit Filters in Countable-Support Symmetric Iterations:\\Corrections and Retained Results}
\author{Frank Gilson}
\date{September 10, 2026}

\begin{document}
		
		\maketitle
		
		\begin{abstract}
			We isolate an algebraic limit-filter construction for countable-support
			symmetric iterations built from coherent successor-stage symmetric systems.
			At limits of uncountable cofinality, countably many generators and their
			conjugating automorphisms are bounded below the limit; at limits of
			cofinality $\omega$, countable-intersection closure is imposed in the
			definition.  The resulting limit filters are normal and
			$\omega_1$-complete.  This gives a ground-coded syntactic closure result:
			the canonical tuple of a ground-model countable sequence of hereditarily
			symmetric names is hereditarily symmetric.  We also recall the standard
			symmetric-extension theorem yielding $\ZF$ from a genuine symmetric system.
			
			The previously claimed preservation of $\DC$ and the unordered-pairs
			application are withdrawn.  Filter completeness alone does not turn the
			ground-coded tuple statement into closure under arbitrary countable sequences
			in the forcing extension.  Moreover, every nonempty finite set of reals has a
			lexicographically least element in $\ZF$, so a family of two-element sets of
			reals without a choice function cannot exist.  The present paper retains only
			the corrected limit-filter observations and the explicitly qualified
			symmetric-name consequences.
		\end{abstract}

	\section*{Correction and current status}
	The principal application previously asserted in this paper is false.  In
	$\ZF$, the lexicographic ordering of $2^\omega$ canonically selects the least
	member of every nonempty finite set of reals.  Consequently, a family of
	two-element sets of reals cannot lack a choice function as was claimed.  In
	addition, under the stated step filter the displayed names for the individual
	reals are not symmetric, so the corresponding pair names are not
	hereditarily symmetric.

	The earlier proof of general $\DC$ preservation also contained a distinct
	gap: it selected, inside $V[G]$, a countable sequence of hereditarily symmetric
	names and then applied a theorem that only concerns countable sequences of
	names already coded in the ground model.  Filter completeness alone does not
	justify that passage.  The general $\DC$ theorem, its consequences, the
	unordered-pairs application, and the finite-support comparison are therefore
	withdrawn.  The corrected algebraic limit-filter results and the standard
	$\ZF$ symmetric-extension facts remain.
	
	\section{Introduction}
	
	Symmetric extensions produce intermediate models of $\ZF$ between a ground
	model $V$ and a forcing extension $V[G]$ by restricting attention to the class
	of hereditarily symmetric names.  Karagila developed a general framework for
	finite-support symmetric iterations \cite{karagila2019}.  The narrower purpose
	of this paper is to record a corrected limit-stage filter construction suited
	to a coherent countable-support presentation.
	
	\paragraph{What is proved here.}
	Working in a fixed background universe $V$, we assume the standard
	successor-step construction of symmetric systems
	(as in \cite[\S\S 3--5]{karagila2019}) and isolate the additional
	\emph{limit-stage} technology needed for \emph{countable-support}
	iterations.  We do not claim here that every proposed iteration template
	automatically forms a symmetric system; that structural verification remains
	part of the hypotheses.
	At each limit $\lambda$ we define a canonical limit filter
	$\tilde{\mathcal F}_\lambda$ extending the head pullbacks of earlier-stage
	filters and prove:
	\begin{itemize}
		\item $\tilde{\mathcal F}_\lambda$ is normal and $\omega_1$--complete (Theorem~\ref{thm:filter-generation});
		\item any ground-model countable sequence of HS names has an HS canonical
		tuple name (Theorem~\ref{thm:hs-closure});
		\item every genuine symmetric system yields the standard transitive
		$\ZF$ symmetric extension (Theorem~\ref{thm:zf-preserved}).
	\end{itemize}
	
	\paragraph{Why the limit filters matter.}
	The central point is that countable-support iterations require
	$\omega_1$-complete symmetry filters at limit stages.
	Its consequence here is syntactic: intersections of the stabilizers of a
	ground-coded countable family remain in the filter.  This observation must not
	be confused with the assertion that the resulting symmetric model is closed
	under all countable sequences appearing in $V[G]$.

	\paragraph{What is \emph{not} claimed.}
	We do not claim here a complete general theory of countable-support
	symmetric iterations comparable in scope to broader support frameworks.
	The successor-step construction is taken from the standard theory, and the
	retained contribution of the present paper is the limit-stage filter
	calculation and its ground-coded syntactic consequence.
	
	We make no dependent-choice preservation claim.  Such conclusions require
	additional hypotheses on the forcing or on how sequences of names are
	captured; see, for example, \cite{karagila2019dc,Banerjee2023}.  We also do
	not address cardinal preservation or other fine structure beyond what is
	explicitly proved.
	All uses of ``countable'' refer to countability in the background universe
	$V$ (Remark~\ref{rem:meta}); for the only countability closure needed here
	(unions of countable supports, i.e.\ subsets of ordinals), see
	Remark~\ref{rem:countable-unions-ord}.
	
	\paragraph{Organization.}
	Section~\ref{sec:prelims-defs} fixes forcing/symmetry notation and records the Symmetry Lemma and related invariance facts used later.
	Section~\ref{sec:count-support-symm} develops the countable-support iteration,
	the limit filters, and the qualified closure statement.  Section~\ref{sec:withdrawn}
	records the withdrawn claims and the obstructions to the earlier application.
	
	\section{Preliminaries: forcing, automorphisms, and symmetric names}\label{sec:prelims-defs}
	
	\subsection{Background theory and metatheoretic conventions}
	\begin{remark}[Metatheory vs.\ object theory]\label{rem:meta}
		All constructions of forcing notions, names, automorphism groups, filters, and
		supports are carried out in a fixed \emph{background} universe $V$.
		
		For the set-length development (existence of the iteration, the limit
		filters, and the standard $\ZF$ symmetric-extension conclusion) we work in
		$V\models\mathsf{ZF}$.  No dependent-choice conclusion is asserted.
		
		When we say a set is \emph{countable}, this is always in the sense of the background universe $V$ (equivalently: there is an injection into $\omega$ in $V$).
		
		No claim is made that these sets are countable inside intermediate symmetric models.
	\end{remark}
	
	\begin{remark}[Countable unions of countable supports in $\ZF$]\label{rem:countable-unions-ord}
		All countability bookkeeping in Sections~\ref{sec:count-support-symm}--\ref{sec:dc-choice} concerns subsets of ordinals (supports and stage indices).
		In $\ZF$, if $\langle A_n:n<\omega\rangle$ is a sequence of countable subsets of an ordinal $\lambda$, then $\bigcup_{n<\omega}A_n$ is countable.
		Indeed, each $A_n\subseteq\lambda$ admits a canonical increasing enumeration (defined by recursion using the well-ordering of $\lambda$, padded by $0$ if $A_n$ is finite),
		and then $(n,m)\mapsto e_n(m)$ is a surjection from $\omega\times\omega$ onto $\bigcup_{n<\omega}A_n$.
		Thus no separate choice principle about countable unions of arbitrary sets is used for the set-length iteration technology here.
	\end{remark}
	
	\subsection{Forcing and names}
	We use standard forcing conventions. A forcing notion is a poset $\mathbb P=(P,\leq)$ with top condition $1_{\mathbb P}$.
	A $\mathbb P$-name is a set $\dot x$ of pairs $\langle \dot y,p\rangle$ with $\dot y$ a $\mathbb P$-name and $p\in\mathbb P$.
	For a $V$-generic filter $G\subseteq\mathbb P$, $\dot x^G$ denotes the usual valuation.
	
	\subsection{Automorphisms and symmetric systems}
	\begin{definition}[Automorphisms]\label{def:autP}
		An \emph{automorphism} of $\mathbb P$ is a bijection $\pi:P\to P$ satisfying
		$p\leq q \iff \pi(p)\leq \pi(q)$ and $\pi(1_{\mathbb P})=1_{\mathbb P}$.
		If $\mathcal G\leq\mathrm{Aut}(\mathbb P)$ is a subgroup, we write $\pi\cdot p$ for $\pi(p)$.
	\end{definition}
	
	\begin{definition}[Normal filter of subgroups]\label{def:normal-filter}
		Let $\mathcal G$ be a group. A \emph{filter of subgroups} $\mathcal F$ on $\mathcal G$ is a nonempty collection of subgroups of $\mathcal G$
		such that:
		\begin{enumerate}
			\item $\mathcal G\in\mathcal F$;
			\item if $H\in\mathcal F$ and $H\leq K\leq\mathcal G$, then $K\in\mathcal F$;
			\item if $H,K\in\mathcal F$, then $H\cap K\in\mathcal F$.
		\end{enumerate}
		It is \emph{normal} if for all $H\in\mathcal F$ and all $\pi\in\mathcal G$, $\pi H\pi^{-1}\in\mathcal F$.
		It is \emph{$\omega_1$-complete} if $\bigcap_{n<\omega}H_n\in\mathcal F$ whenever $H_n\in\mathcal F$ for all $n<\omega$.
	\end{definition}
	
	\begin{definition}[$\omega_1$-completion of a normal filter]\label{def:omega1-completion}
		Let $\mathcal F$ be a normal filter of subgroups on $\mathcal G$.
		Define its \emph{$\omega_1$-completion} $\mathcal F^{\omega_1}$ by
		\[
		K\in\mathcal F^{\omega_1}\quad\Longleftrightarrow\quad
		\exists\langle H_n:n<\omega\rangle\in\mathcal F^\omega\ \Bigl(\ \bigcap_{n<\omega}H_n\le K\ \Bigr).
		\]
	\end{definition}
	
	\begin{lemma}[$\omega_1$-completion is the smallest normal $\omega_1$-complete extension]\label{lem:omega1-completion}
		If $\mathcal F$ is a normal filter of subgroups on $\mathcal G$, then $\mathcal F^{\omega_1}$ is a normal, $\omega_1$-complete filter of subgroups on $\mathcal G$,
		and it is the smallest normal, $\omega_1$-complete filter extending $\mathcal F$.
	\end{lemma}
	
	\begin{proof}
		Upward closure and containing $\mathcal G$ are immediate.
		$\omega_1$-completeness holds because if $K_m\in\mathcal F^{\omega_1}$ is witnessed by $\langle H^m_n:n<\omega\rangle\subseteq\mathcal F$,
		then $\bigcap_m K_m$ is witnessed by the single $\mathcal F$-sequence obtained by concatenating the sequences $\langle H^m_n:n<\omega\rangle$ over $m<\omega$.
		Normality is preserved since for $\pi\in\mathcal G$,
		\[
		\pi\Bigl(\bigcap_{n<\omega}H_n\Bigr)\pi^{-1}=\bigcap_{n<\omega}\pi H_n\pi^{-1},
		\]
		and $\mathcal F$ is normal.
		Minimality is immediate: any normal $\omega_1$-complete filter extending $\mathcal F$ must contain every countable intersection of members of $\mathcal F$,
		hence must contain $\mathcal F^{\omega_1}$.
	\end{proof}
	
	\begin{definition}[Symmetric system]\label{def:symm-system}
		A \emph{symmetric system} is a triple $\langle\mathbb P,\mathcal G,\mathcal F\rangle$ where
		$\mathbb P$ is a forcing notion, $\mathcal G\leq\mathrm{Aut}(\mathbb P)$, and $\mathcal F$ is a normal filter of subgroups of $\mathcal G$.
	\end{definition}
	
	\subsection{Action on names and the Symmetry Lemma}
	Given $\langle\mathbb P,\mathcal G,\mathcal F\rangle$, define recursively the action of $\pi\in\mathcal G$ on $\mathbb P$-names by
	\[
	\pi\dot x\ :=\ \{\langle \pi\dot y,\ \pi p\rangle : \langle \dot y,p\rangle\in\dot x\}.
	\]
	
	\begin{definition}[Stabilizers]\label{def:stabilizers}
		For a name $\dot x$, its \emph{(setwise) stabilizer} is
		\[
		\sym_{\mathcal G}(\dot x)\ :=\ \{\pi\in\mathcal G : \pi\dot x=\dot x\}.
		\]
		When $\mathcal G$ is clear from context, we write $\sym(\dot x)$.
	\end{definition}
	
	\begin{lemma}[Symmetry Lemma / Permutation Lemma]\label{lem:symmetry-lemma}
		For any formula $\varphi$ in the forcing language, any $\pi\in\mathcal G$, any $p\in\mathbb P$,
		and any names $\dot x_1,\dots,\dot x_n$,
		\[
		p\Vdash_{\mathbb P}\varphi(\dot x_1,\dots,\dot x_n)\quad\Longleftrightarrow\quad
		\pi p\Vdash_{\mathbb P}\varphi(\pi\dot x_1,\dots,\pi\dot x_n).
		\]
		See, e.g., \cite[Lemma~2.4]{karagila2019}.
	\end{lemma}
	
	\subsection{Symmetric and hereditarily symmetric names}
	\begin{definition}[Symmetric / hereditarily symmetric]\label{def:HS}
		Let $\langle\mathbb P,\mathcal G,\mathcal F\rangle$ be a symmetric system.
		A name $\dot x$ is \emph{$\mathcal F$-symmetric} if $\sym_{\mathcal G}(\dot x)\in\mathcal F$.
		A name $\dot x$ is \emph{hereditarily $\mathcal F$-symmetric} if $\dot x$ is $\mathcal F$-symmetric and
		every $\dot y$ appearing in the transitive closure of $\dot x$ (equivalently, every $\dot y\in\mathrm{dom}(\dot x)$ and recursively)
		is hereditarily $\mathcal F$-symmetric.
		We write $\HS(\mathbb P,\mathcal G,\mathcal F)$ for the class of hereditarily $\mathcal F$-symmetric $\mathbb P$-names.
	\end{definition}
	
	\subsection{Dependent Choice (set formulation)}
	\begin{definition}[Dependent Choice]\label{def:DC}
		Let $\mu$ be a nonzero ordinal. $\DC_\mu$ is the statement:
		for every \emph{set} $A$ and every binary relation $R\subseteq A\times A$ such that
		$\forall x\in A\,\exists y\in A\ (xRy)$, there exists a sequence $\langle a_\xi:\xi<\mu\rangle$ of elements of $A$
		with $a_\xi R a_{\xi+1}$ for all $\xi+1<\mu$.
		We write $\DC=\DC_\omega$ and $\DC_{<\kappa}$ for $\forall\mu<\kappa\,\DC_\mu$.
	\end{definition}
	
	\subsection{Two-step iterations (factorization)}
	\begin{fact}[Two-step factorization]\label{fact:two-step}
		Let $\mathbb P$ be a forcing notion and let $\dot{\mathbb Q}$ be a $\mathbb P$-name for a forcing notion.
		Let $\mathbb P*\dot{\mathbb Q}$ be the standard two-step iteration. If $G*H\subseteq \mathbb P*\dot{\mathbb Q}$ is $V$-generic,
		then $G\subseteq\mathbb P$ is $V$-generic and $H\subseteq \dot{\mathbb Q}^G$ is $V[G]$-generic, and conversely every such pair arises.
	\end{fact}
	
	\begin{remark}[Notation discipline for later sections]\label{rem:notation-discipline}
		In the iteration sections we will \emph{not} overload projection symbols.
		For forcing, we use restriction $p\upharpoonright\beta$ and (when needed) a map $\pi_{\beta,\lambda}:\mathbb P_\lambda\to\mathbb P_\beta$.
		For groups, we use a distinct restriction homomorphism $\rho_{\beta,\lambda}:\mathcal G_\lambda\to\mathcal G_\beta$.
	\end{remark}
	
	\section{Limit-Stage Filters for Countable-Support Symmetric Iterations}\label{sec:count-support-symm}
	
	\subsection{Framework and Setup}
	
	\begin{definition}[Countable-support iteration setup]\label{def:cs-iteration}
		Fix an ordinal $\Theta$.
		A \emph{countable-support symmetric iteration of length $\Theta$} is a system
		\[
		\langle (\mathbb P_\alpha,\mathcal G_\alpha,\mathcal F_\alpha):\alpha\le\Theta\rangle
		\]
		together with, for each $\alpha<\Theta$, a $\mathbb P_\alpha$-name
		\[
		\dot{\mathscr S}_\alpha=\langle \dot{\mathbb Q}_\alpha,\dot{\mathcal H}_\alpha,\dot{\mathcal K}_\alpha\rangle
		\]
		such that
		\[
		1_{\mathbb P_\alpha}\Vdash_{\mathbb P_\alpha}\text{``$\dot{\mathscr S}_\alpha$ is a symmetric system.''}
		\]
		The recursion is as follows.
		
		\paragraph{Stage $0$.}
		$\mathbb P_0=\{1\}$, $\mathcal G_0=\{\mathrm{id}\}$, and $\mathcal F_0=\{\mathcal G_0\}$.
		
		\paragraph{Successor stage $\alpha+1$.}
		Let $\dot{\mathscr S}_\alpha=\langle \dot{\mathbb Q}_\alpha,\dot{\mathcal H}_\alpha,\dot{\mathcal K}_\alpha\rangle$ be as above.
		Define
		\[
		\mathbb P_{\alpha+1}=\mathbb P_\alpha*\dot{\mathbb Q}_\alpha.
		\]
		Define $\mathcal G_{\alpha+1}$ and an intermediate \emph{normal} successor-stage filter
		$\mathcal F^{\mathrm{Kar}}_{\alpha+1}$ by the \emph{standard two-step symmetric-system construction}
		from $(\mathbb P_\alpha,\mathcal G_\alpha,\mathcal F_\alpha)$ and the name $\dot{\mathscr S}_\alpha$:
		there is a canonical restriction homomorphism (``head restriction'')
		\[
		\rho_{\alpha,\alpha+1}:\mathcal G_{\alpha+1}\to\mathcal G_\alpha
		\]
		and $\mathcal F^{\mathrm{Kar}}_{\alpha+1}$ is a normal filter of subgroups of $\mathcal G_{\alpha+1}$ extending the pullbacks
		$\rho_{\alpha,\alpha+1}^{-1}(H)$ for $H\in\mathcal F_\alpha$ and the tail symmetries coming from $\dot{\mathcal K}_\alpha$.
		(We do not reprove this successor-step construction here; see, e.g., \cite[\S\S 3--5]{karagila2019}.)
		
		For the countable-support limit-stage arguments below, we require $\omega_1$--completeness at intermediate stages, so we define
		\[
		\mathcal F_{\alpha+1}\ :=\ \bigl(\mathcal F^{\mathrm{Kar}}_{\alpha+1}\bigr)^{\omega_1},
		\]
		the $\omega_1$--completion of $\mathcal F^{\mathrm{Kar}}_{\alpha+1}$ in the sense of Definition~\ref{def:omega1-completion}
		(which is the smallest normal $\omega_1$--complete extension by Lemma~\ref{lem:omega1-completion}).
		
		\paragraph{Limit stage $\lambda\le\Theta$ with $\cf(\lambda)=\omega$.}
		$\mathbb P_\lambda$ is the usual countable-support limit of $\langle \mathbb P_\alpha,\dot{\mathbb Q}_\alpha:\alpha<\lambda\rangle$:
		conditions are partial functions $p$ with countable domain $\dom(p)\subseteq\lambda$, where for each $\beta\in\dom(p)$,
		$p(\beta)$ is a $\mathbb P_\beta$-name for a condition in $\dot{\mathbb Q}_\beta$, ordered by
		$q\le p$ iff $\dom(q)\supseteq\dom(p)$ and for all $\beta\in\dom(p)$,
		\[
		q\upharpoonright\beta\ \Vdash_{\mathbb P_\beta}\ q(\beta)\le p(\beta).
		\]
		Define $\supp(p):=\dom(p)$.
		
		$\mathcal G_\lambda$ is the inverse limit of $\langle\mathcal G_\beta,\rho_{\gamma,\beta}:\gamma<\beta<\lambda\rangle$ and acts
		coordinatewise on $\mathbb P_\lambda$ (Lemma~\ref{lem:limit-action}).
		The limit filter $\tilde{\mathcal F}_\lambda$ is defined in Definition~\ref{def:limit-filter}.
		
		\paragraph{Limit stage $\lambda\le\Theta$ with $\cf(\lambda)\ge\omega_1$.}
		We identify the countable-support limit with the direct limit (Remark~\ref{rem:limit-bookkeeping}):
		\[
		\mathbb P_\lambda=\bigcup_{\beta<\lambda}\mathbb P_\beta,\qquad
		\mathcal G_\lambda=\bigcup_{\beta<\lambda}\mathcal G_\beta,
		\]
		with the canonical restriction maps $\pi_{\beta,\lambda}(p)=p\upharpoonright\beta$ and $\rho_{\beta,\lambda}$.
		The limit filter $\tilde{\mathcal F}_\lambda$ is the normal filter generated by head pullbacks from earlier stages
		(Definition~\ref{def:limit-filter}), and it is $\omega_1$-complete by stage-bounding
		(Remark~\ref{rem:dl-omega1-completeness}).
	\end{definition}
	
	\begin{remark}[Bookkeeping at uncountable cofinality]\label{rem:limit-bookkeeping}
		Let $\lambda$ be a limit ordinal with $\cf(\lambda)\ge\omega_1$, and let $\mathbb P_\lambda$ be the usual countable-support
		limit of $\langle \mathbb P_\alpha,\dot{\mathbb Q}_\alpha:\alpha<\lambda\rangle$.
		Then every $p\in\mathbb P_\lambda$ has bounded support: since $\supp(p)$ is countable and $\cf(\lambda)>\omega$, there is
		$\beta<\lambda$ with $\supp(p)\subseteq\beta$, hence $p\in\mathbb P_\beta$.
		Therefore $\mathbb P_\lambda=\bigcup_{\beta<\lambda}\mathbb P_\beta$ (direct limit identification),
		and similarly $\mathcal G_\lambda=\bigcup_{\beta<\lambda}\mathcal G_\beta$ with the evident restriction maps.
	\end{remark}
	
	\begin{lemma}[Limit-stage coordinatewise action]\label{lem:limit-action}
		Let $\lambda$ be a limit stage with $\cf(\lambda)=\omega$.
		If $g=\langle g_\beta:\beta<\lambda\rangle\in\mathcal G_\lambda$, define $(g\cdot p)(\beta)=g_\beta\cdot p(\beta)$ for $\beta\in\dom(p)$.
		Then $p\mapsto g\cdot p$ is an automorphism of $\mathbb P_\lambda$.
	\end{lemma}
	
	\begin{proof}
		It is immediate that $\dom(g\cdot p)=\dom(p)$ and $g\cdot 1_{\mathbb P_\lambda}=1_{\mathbb P_\lambda}$.
		For order preservation, suppose $q\le p$ in $\mathbb P_\lambda$.
		Fix $\beta\in\dom(p)$. Then $q\upharpoonright\beta\Vdash_{\mathbb P_\beta} q(\beta)\le p(\beta)$.
		Apply the invariance of the forcing relation under automorphisms from Section~\ref{sec:prelims-defs}
		(Lemma~\ref{lem:symmetry-lemma} in its forcing-theorem form): since $g\upharpoonright\beta\in\mathcal G_\beta$,
		\[
		(g\cdot q)\upharpoonright\beta\ \Vdash_{\mathbb P_\beta}\ (g\cdot q)(\beta)\le (g\cdot p)(\beta).
		\]
		As this holds for all $\beta\in\dom(p)$, we have $g\cdot q\le g\cdot p$.
		Bijectivity follows from applying the same argument to $g^{-1}$.
	\end{proof}
	
	\begin{remark}[Comparison with related frameworks]
		\label{rem:karagila-framework-relation}
		The successor-step symmetric-system construction used above is standard 
		(see, e.g., \cite[\S\S 3--5]{karagila2019}).
		What is new in the present paper is the \emph{limit-stage} treatment needed 
		for countable support: at $\cf(\lambda)\ge\omega_1$ the direct-limit 
		identification (Remark~\ref{rem:limit-bookkeeping}) allows stage-bounding, 
		while at $\cf(\lambda)=\omega$ we define the limit filter by taking the smallest 
		normal $\omega_1$-complete extension of the head pullbacks 
		(Definition~\ref{def:limit-filter}).
		
		The independent concurrent work of Karagila--Schilhan \cite{KaragilaSchilhan2026} 
		develops a broad structural theory of symmetric extensions---including products, 
		reduced products, quotients, and a general $I$-support iteration framework 
		(where $I$ is any ideal on the index ordinal)---and proves results on 
		Kinna--Wagner principles and the HOD theorem.
		Their $I$-support definition subsumes countable support as a special case.
		The present paper is complementary in focus: it records the specific
		$\omega_1$-complete normal limit-filter calculation and its syntactic
		consequence.  It should be read as a limit-stage supplement to
		the standard successor-step theory, rather than as a replacement for more
		general support frameworks.
	\end{remark}
	
	\begin{remark}[Uncountable cofinality: stage-bounding yields \texorpdfstring{$\omega_1$}{omega1}-completeness]
		\label{rem:dl-omega1-completeness}
		Assume $\cf(\lambda)\ge\omega_1$ and work in the direct-limit identification of
		Remark~\ref{rem:limit-bookkeeping}.
		Let $\tilde{\mathcal F}_\lambda$ be the limit filter from Definition~\ref{def:limit-filter}(1),
		i.e.\ the smallest normal filter on $\mathcal G_\lambda$ containing the head pullbacks
		$\rho_{\beta,\lambda}^{-1}(H)$ for $\beta<\lambda$ and $H\in\mathcal F_\beta$.
		
		Then $\tilde{\mathcal F}_\lambda$ is $\omega_1$-complete.
		Indeed, let $\langle K_n:n<\omega\rangle\subseteq\tilde{\mathcal F}_\lambda$.
		For each $n$, by the normal-generation clause in Definition~\ref{def:limit-filter}(1),
		fix witnesses
		$k_n<\omega$, ordinals $\beta_{n,i}<\lambda$, groups $H_{n,i}\in\mathcal F_{\beta_{n,i}}$,
		and elements $g_{n,i}\in\mathcal G_\lambda$ such that
		\[
		\bigcap_{i<k_n} g_{n,i}\,\rho_{\beta_{n,i},\lambda}^{-1}(H_{n,i})\,g_{n,i}^{-1}\ \le\ K_n.
		\]
		By the direct-limit hypothesis, for each $n,i$ there is
		$\gamma_{n,i}<\lambda$ such that $g_{n,i}$ is the canonical extension of an
		element of $\mathcal G_{\gamma_{n,i}}$.  Since
		$\cf(\lambda)\ge\omega_1$, the countable set of all ordinals
		$\beta_{n,i}$ and $\gamma_{n,i}$ is bounded in $\lambda$; let
		$\beta^*<\lambda$ bound them all.
		
		For each $n,i$, set $H'_{n,i}:=\rho_{\beta_{n,i},\beta^*}^{-1}(H_{n,i})$.
		Then $H'_{n,i}\in\mathcal F_{\beta^*}$ by Lemma~\ref{lem:push-restrict-preserve}, and
		\[
		\rho_{\beta^*,\lambda}^{-1}(H'_{n,i})\ \le\ \rho_{\beta_{n,i},\lambda}^{-1}(H_{n,i})
		\]
		by coherence of restrictions.
		Replacing each head pullback by $\rho_{\beta^*,\lambda}^{-1}(H'_{n,i})$ only strengthens the left-hand side, so
		we may assume all generators come from stage $\beta^*$.
		
		Each $g_{n,i}$ is now represented by some
		$h_{n,i}\in\mathcal G_{\beta^*}$.  Normality puts
		$h_{n,i}H'_{n,i}h_{n,i}^{-1}$ in $\mathcal F_{\beta^*}$.
		Use $\omega_1$-completeness of $\mathcal F_{\beta^*}$ to intersect these
		countably many conjugated stage-$\beta^*$ subgroups, and pull the result
		back along $\rho_{\beta^*,\lambda}$.
		This produces an element of $\tilde{\mathcal F}_\lambda$ contained in $\bigcap_{n<\omega}K_n$.
	\end{remark}
	
	\subsection{Filter Generation at Limit Stages}
	\label{ss:filter-generation}
	
	Let $\lambda$ be a limit ordinal and assume inductively that each $\mathcal F_\beta$ for $\beta<\lambda$ is a normal, $\omega_1$-complete filter on $\mathcal G_\beta$.
	We define the limit filter $\tilde{\mathcal F}_\lambda$ uniformly from the head pullbacks $\rho_{\beta,\lambda}^{-1}(H)$ (Definition~\ref{def:limit-filter}).
	At $\cf(\lambda)\ge\omega_1$ we take the smallest normal filter containing these generators, and $\omega_1$-completeness follows by stage-bounding (Remark~\ref{rem:dl-omega1-completeness}).
	At $\cf(\lambda)=\omega$ we instead take the smallest normal $\omega_1$-complete filter extending the same generator family.
	
	\begin{remark}[Head pullbacks of subgroups]\label{rem:head-pullback}
		For $\beta<\lambda$ and $H\le\mathcal G_\beta$, we write
		\[
		\rho_{\beta,\lambda}^{-1}(H)=\{g\in\mathcal G_\lambda:\rho_{\beta,\lambda}(g)\in H\}.
		\]
		When pulling back \emph{filters} along homomorphisms we use the notation $\pi^*\mathcal F$ from
		Lemma~\ref{lem:push-restrict-preserve}.
	\end{remark}
	
	\begin{lemma}[Pullbacks preserve normality and \texorpdfstring{$\omega_1$}{omega1}-completeness]
		\label{lem:push-restrict-preserve}
		Let $\pi:\mathcal G\to\mathcal H$ be a homomorphism and let $\mathcal F$ be a normal, $\omega_1$-complete filter on $\mathcal H$.
		Define the pullback filter on $\mathcal G$ by
		\[
		\pi^{*}\mathcal F\ :=\ \{\,K\le \mathcal G:\ \exists H\in\mathcal F\ \text{with}\ \pi^{-1}(H)\le K\,\}.
		\]
		Then $\pi^{*}\mathcal F$ is a normal, $\omega_1$-complete filter on $\mathcal G$.
	\end{lemma}
	
	\begin{proof}
		\emph{Pullback.}
		Upward closure and containing $\mathcal G$ are immediate.
		For $\omega_1$-completeness, suppose $\langle K_n:n<\omega\rangle$ with each $K_n\in\pi^{*}\mathcal F$.
		Choose $H_n\in\mathcal F$ with $\pi^{-1}(H_n)\le K_n$.
		Then $\bigcap_{n<\omega}H_n\in\mathcal F$ by $\omega_1$-completeness, and
		\[
		\pi^{-1}\!\left(\bigcap_{n<\omega}H_n\right)=\bigcap_{n<\omega}\pi^{-1}(H_n)\ \le\ \bigcap_{n<\omega}K_n,
		\]
		so $\bigcap_{n<\omega}K_n\in\pi^{*}\mathcal F$.
		For normality: if $K\in\pi^{*}\mathcal F$ witnessed by $H\in\mathcal F$ and $g\in\mathcal G$, then
		$\pi(g)H\pi(g)^{-1}\in\mathcal F$ and
		\[
		\pi^{-1}(\pi(g)H\pi(g)^{-1}) = g\,\pi^{-1}(H)\,g^{-1}\ \le\ gKg^{-1},
		\]
		so $gKg^{-1}\in\pi^{*}\mathcal F$.
		This proves the assertion.
	\end{proof}
	
	\begin{definition}[Supports of conditions and names (countable-support case)]
		\label{def:support-recursive}
		Let $\lambda$ be a limit stage.
		
		\smallskip
		\noindent\emph{Support of a condition.}
		For $p\in\mathbb P_\lambda$, let $\supp(p)\subseteq\lambda$ be the set of coordinates at which $p$ is nontrivial
		(in the usual sense for the iteration/product presentation at $\lambda$).
		In the countable-support presentation, $\supp(p)$ is required to be countable.
		
		\smallskip
		\noindent\emph{Support of a name.}
		Define the support of a $\mathbb P_\lambda$--name $\dot x$ by recursion on name-rank:
		\[
		\supp(\dot x)\ :=\ \bigcup\Bigl\{\ \supp(p)\ \cup\ \supp(\sigma)\ :\ (\sigma,p)\in\dot x\ \Bigr\}.
		\]
		We say that $\dot x$ has \emph{countable support} if $\supp(\dot x)$ is countable.
	\end{definition}
	
	\begin{definition}[Limit filter]\label{def:limit-filter}
		Let $\lambda$ be a limit ordinal and assume we have defined $(\mathbb P_\beta,\mathcal G_\beta,\mathcal F_\beta)$ for all $\beta<\lambda$,
		with coherent restriction homomorphisms $\rho_{\beta,\lambda}:\mathcal G_\lambda\to\mathcal G_\beta$.
		
		Define the \emph{head-generator family} at $\lambda$ by
		\[
		\mathcal B_\lambda\ :=\ \bigl\{\ \rho_{\beta,\lambda}^{-1}(H)\ :\ \beta<\lambda\ \&\ H\in\mathcal F_\beta\ \bigr\}.
		\]
		
		\begin{enumerate}
			\item If $\cf(\lambda)\ge\omega_1$, let $\tilde{\mathcal F}_\lambda$ be the normal filter on $\mathcal G_\lambda$ \emph{generated by} $\mathcal B_\lambda$, equivalently:
			\[
			K\in\tilde{\mathcal F}_\lambda\ \Longleftrightarrow\
			\exists k<\omega\ \exists B_0,\dots,B_{k-1}\in\mathcal B_\lambda\ \exists g_0,\dots,g_{k-1}\in\mathcal G_\lambda\
			\Bigl(\bigcap_{i<k} g_iB_i g_i^{-1}\le K\Bigr).
			\]
			
			\item If $\cf(\lambda)=\omega$, let $\tilde{\mathcal F}_\lambda$ be the normal $\omega_1$--complete filter on $\mathcal G_\lambda$ \emph{generated by} $\mathcal B_\lambda$, equivalently:
			\[
			K\in\tilde{\mathcal F}_\lambda\ \Longleftrightarrow\
			\exists\langle B_n:n<\omega\rangle\in\mathcal B_\lambda^\omega\ \exists\langle g_n:n<\omega\rangle\in\mathcal G_\lambda^\omega\
			\Bigl(\bigcap_{n<\omega} g_nB_n g_n^{-1}\le K\Bigr).
			\]
		\end{enumerate}
	\end{definition}
	
	\begin{remark}[Syntactic support vs.\ background countability]
		By Definition~\ref{def:support-recursive}, $\supp(\dot x)$ is computed recursively from the supports of the conditions
		appearing in $\dot x$ (and from the supports of subnames). When we say “$\dot x$ has countable support’’ in the
		countable-support part of the paper, we mean that $\supp(\dot x)$ is countable in the ambient background universe
		(cf.\ Remark~\ref{rem:meta}).
	\end{remark}
	
	The $\omega_1$--completeness argument splits by cofinality.
	If $\cf(\lambda)\ge\omega_1$, stage-bounding reduces countable intersections in $\tilde{\mathcal F}_\lambda$
	to a single stage filter (Remark~\ref{rem:dl-omega1-completeness}).
	If $\cf(\lambda)=\omega$, $\omega_1$--completeness is built into the definition of $\tilde{\mathcal F}_\lambda$
	as the smallest normal $\omega_1$--complete filter containing the head generators (Definition~\ref{def:limit-filter}).
	
	\begin{lemma}[Stage-Bounding Lemma---uncountable cofinality only]
		\label{lemma:stage-bounding}
		If $\lambda$ is a limit ordinal with $\cf(\lambda)\ge\omega_1$ and $\{\beta_n:n<\omega\}\subset\lambda$ is any countable set of ordinals below $\lambda$, then 
		\[
		\sup_{n<\omega}\beta_n \;<\;\lambda.
		\]
	\end{lemma}
	
	\begin{proof}
		If $\sup_n\beta_n=\lambda$, then $\{\beta_n:n<\omega\}$ would be a countable cofinal subset of $\lambda$, contradicting $\cf(\lambda)\ge\omega_1$. Hence $\sup_n\beta_n<\lambda$.
	\end{proof}
	
	\begin{remark}[Failure of stage-bounding at countable cofinality]
		\label{rem:stage-bounding-fails}
		When $\cf(\lambda)=\omega$, the stage-bounding lemma fails. For example, if $\lambda=\omega$, the set $\{n:n<\omega\}$ is countable but $\sup_n n=\omega=\lambda$. Similarly, if $\lambda=\aleph_\omega$, the sequence $\{\aleph_n:n<\omega\}$ is countable and cofinal. This is precisely why, at limits of countable cofinality, $\omega_1$--completeness must be built into the limit-filter definition (Definition~\ref{def:limit-filter}(2)), rather than proved by stage-bounding.
	\end{remark}
	
	\begin{theorem}[Limit filter properties]\label{thm:filter-generation}
		Let $\lambda$ be a limit stage and let $\tilde{\mathcal F}_\lambda$ be as in Definition~\ref{def:limit-filter}.
		Assume inductively that each $\mathcal F_\beta$ for $\beta<\lambda$ is a normal $\omega_1$--complete filter on $\mathcal G_\beta$.
		Then:
		\begin{enumerate}
			\item $\tilde{\mathcal F}_\lambda$ is a normal filter on $\mathcal G_\lambda$ and contains every head pullback
			$\rho_{\beta,\lambda}^{-1}(H)$ for $\beta<\lambda$ and $H\in\mathcal F_\beta$.
			\item $\tilde{\mathcal F}_\lambda$ is $\omega_1$--complete.
		\end{enumerate}
	\end{theorem}
	
	\begin{proof}
		By Definition~\ref{def:limit-filter}, $\tilde{\mathcal F}_\lambda$ is (by construction) the smallest normal filter (or normal $\omega_1$--complete filter)
		containing the generator family $\mathcal B_\lambda$, so it is a normal filter and contains each $\rho_{\beta,\lambda}^{-1}(H)$.
		
		For $\omega_1$--completeness:
		if $\cf(\lambda)=\omega$, closure under countable intersections is built into the defining criterion in
		Definition~\ref{def:limit-filter}(2) by interleaving witnesses.
		If $\cf(\lambda)\ge\omega_1$, $\omega_1$--completeness follows from stage-bounding as proved in
		Remark~\ref{rem:dl-omega1-completeness}.
	\end{proof}
	
	\subsection{Hereditarily Symmetric Names in Countable Support Iterations}
	\label{sec:HS-countable}
	
	Having defined the limit filter $\tilde{\mathcal F}_\lambda$
	(Definition~\ref{def:limit-filter}) and verified its normality and
	$\omega_1$--completeness (Theorem~\ref{thm:filter-generation}), we isolate
	the syntactic closure statement that follows directly from completeness.
	
	\begin{definition}[HS at stage $\lambda$]\label{def:HS-stage}
		Let $\langle \mathbb P_\lambda,\mathcal G_\lambda,\widetilde{\mathcal F}_\lambda\rangle$ be the symmetric system at stage $\lambda$.
		We write
		\[
		\HS_\lambda := \HS(\mathbb P_\lambda,\mathcal G_\lambda,\widetilde{\mathcal F}_\lambda)
		\]
		for the class of hereditarily $\widetilde{\mathcal F}_\lambda$-symmetric $\mathbb P_\lambda$-names, as in Definition~\ref{def:HS}.
		When $\lambda$ is fixed, we simply write $\HS$.
	\end{definition}
	
	\begin{equation}\label{eq:tuple-sym}
		\sym\!\bigl(\langle \dot x_j : j\in J\rangle\bigr)
		\;=\; \bigcap_{j\in J}\sym(\dot x_j).
	\end{equation}
	% Equality is of subgroups of $\mathcal G_\lambda$. For finite $J$ this yields
	% closure under pairing/products; for $|J|<\kappa$ we combine \eqref{eq:tuple-sym}
	% with $\kappa$-completeness of $\tilde{\mathcal F}_\lambda$ to keep the
	% intersection inside the filter.
	\begin{proof}[Justification of \eqref{eq:tuple-sym}]
		Let $\dot s=\langle \dot x_j:j\in J\rangle$ be the canonical sequence name.
		Then $g\in\sym(\dot s)$ iff $g(\dot s)=\dot s$, which (by the recursive definition of the action on names)
		holds iff $g(\dot x_j)=\dot x_j$ for all $j\in J$, i.e.\ iff $g\in\bigcap_{j\in J}\sym(\dot x_j)$.
	\end{proof}
	
	\begin{theorem}[Ground-coded tuple closure]
		\label{thm:hs-closure}
		Let $\langle\dot x_n:n<\omega\rangle$ be a sequence that belongs to the
		ground model $V$, and suppose every $\dot x_n$ is an HS name.  Then the
		canonical name $\langle\dot x_n:n<\omega\rangle$ for the ordered tuple is
		HS.  In particular, HS names are closed under ground-coded finite tuples.
	\end{theorem}
	
	\begin{proof}
		Let $\dot s$ be the canonical tuple name formed in $V$.  By
		Equation~\eqref{eq:tuple-sym},
		\[
		\sym(\dot s)\ =\ \bigcap_{n<\omega}\sym(\dot x_n).
		\]
		Since each $\sym(\dot x_n)\in\tilde{\mathcal F}_\lambda$ and $\tilde{\mathcal F}_\lambda$ is $\omega_1$--complete
		(Theorem~\ref{thm:filter-generation}), we have
		$\sym(\dot s)\in\tilde{\mathcal F}_\lambda$.  Every constituent of
		$\dot s$ is hereditarily symmetric by hypothesis, so $\dot s$ is HS.
	\end{proof}
	
	\begin{remark}[What the tuple theorem does not say]
		The quantification over sequences in Theorem~\ref{thm:hs-closure} is in
		the ground model.  The theorem does not imply
		$M^\omega\cap V[G]\subseteq M$: an arbitrary sequence of elements of $M$
		chosen in $V[G]$ need not admit a sequence of HS names coded in $V$.
		This distinction is exactly why the withdrawn $\DC$ argument does not follow
		from $\omega_1$-completeness alone.
	\end{remark}
	
	% ----------------------------
	
	\subsection{\texorpdfstring{$\ZF$}{ZF} Preservation}
	\label{sec:zf-preservation}
	
	\begin{theorem}[Standard symmetric-extension theorem]
		\label{thm:zf-preserved}
		Let $V\models\ZF$, let
		$\langle\mathbb P,\mathcal G,\mathcal F\rangle$ be a symmetric system,
		and let $G\subseteq\mathbb P$ be $V$--generic.  If $\HS$ is the class of
		hereditarily $\mathcal F$--symmetric names, then
		\[
		M=\{\dot x^G:\dot x\in\HS\}
		\]
		is a transitive model of $\ZF$ with $V\subseteq M\subseteq V[G]$.
	\end{theorem}

	\begin{proof}
		This is the standard Symmetric Extension Theorem; see, for example,
		\cite[Chapter~15]{Jech03} or the modern iteration treatment
		\cite{karagila2019}.  Its proof uses the relativized forcing relation for
		hereditarily symmetric names.  In particular, Separation and Replacement
		are not obtained merely by restricting the ordinary forcing relation to
		$M$.  Normality of $\mathcal F$ is the relevant filter hypothesis for the
		standard $\ZF$ theorem; $\omega_1$-completeness is not required for this
		conclusion.
	\end{proof}

	\begin{remark}[Scope]
		Theorem~\ref{thm:zf-preserved} applies once the forcing, automorphism group,
		and normal filter actually form a symmetric system.  The algebraic
		construction of Theorem~\ref{thm:filter-generation} supplies the filter
		properties at a limit stage under its stated coherence hypotheses, but it
		does not replace the remaining verification that a proposed iteration
		template is a genuine symmetric system.
	\end{remark}

	% ----------------------------
	
	\subsection{Withdrawn dependent-choice claim}
	\label{sec:dc-choice}

	\paragraph{Withdrawn claim.}
	The assertion that $\omega_1$-completeness of the symmetry filter alone
	implies $\DC$ in the symmetric model is withdrawn.  The earlier argument
	worked in $V[G]$, obtained a sequence $\langle a_n:n<\omega\rangle$ of
	elements of $M$, and chose HS names $\dot a_n$ for its entries.  It then
	applied Theorem~\ref{thm:hs-closure}.  This is invalid because the sequence
	$\langle\dot a_n:n<\omega\rangle$ selected in $V[G]$ need not be a
	ground-model sequence of names.  Theorem~\ref{thm:hs-closure} is a
	ground-coded syntactic closure theorem, not a proof that
	$M^\omega\cap V[G]\subseteq M$.

	Consequently, no conclusion about $\DC$, Dedekind-finite sets, or amorphous
	sets is drawn here from filter completeness alone.  Preservation theorems for
	dependent choice require additional hypotheses, commonly closure or chain
	conditions on the forcing together with the appropriate completeness of the
	filter; see \cite{karagila2019dc,Banerjee2023}.

	\subsection{Iteration length}\label{subsec:length}
	
	\begin{proposition}[Set-length countable-support iterations]\label{prop:set-length}
		Fix a countable-support iteration template as in Section~\ref{sec:count-support-symm}.
		Then for every ordinal $\Theta$ the recursion defining
		\[
		\langle (\mathbb P_\alpha,\mathcal G_\alpha,\mathcal F_\alpha):\alpha\le\Theta\rangle
		\]
		is well-defined. For each limit $\lambda\le\Theta$, the induced limit filter $\tilde{\mathcal F}_\lambda$
		is normal and $\omega_1$-complete.  Whenever the stage data form a genuine
		symmetric system, the corresponding symmetric model $M_\lambda$ satisfies
		$\ZF$ by the standard symmetric-extension theorem.
	\end{proposition}
	
	\begin{proof}
		The successor step is the standard two-step symmetric-system construction (as cited in Definition~\ref{def:cs-iteration}).
		At limit stages $\lambda\le\Theta$, we use Remark~\ref{rem:limit-bookkeeping} to identify the countable-support limit with the direct limit
		when $\cf(\lambda)\ge\omega_1$, and we use the usual countable-support definition when $\cf(\lambda)=\omega$.
		In both cofinality cases, the limit filter is defined by Definition~\ref{def:limit-filter} and is normal and $\omega_1$-complete
		by Theorem~\ref{thm:filter-generation}.
		The symmetric model at stage $\lambda$ satisfies $\ZF$ by
		Theorem~\ref{thm:zf-preserved}, provided the construction supplies the
		required symmetric system.  No dependent-choice conclusion is asserted.
	\end{proof}
	
	\begin{remark}[Class-length iterations]
		\label{rem:class-length}
		We do not claim a general class-length $\ZF$-preservation theorem
		in this paper.
		The set-length results (Proposition~\ref{prop:set-length} and
		Theorem~\ref{thm:zf-preserved})
		are stated and proved for set-length iterations only, where the
		ambient forcing extension is set-forcing based.
		
		At class length, $\ZF$ preservation for symmetric submodels requires
		hypotheses beyond pretameness: there are class forcings that are
		pretame yet whose symmetric submodels fail Collection or Power Set.
		Identifying the correct general hypothesis is a delicate open problem
		in class forcing theory.
		
		Specific class-length symmetric iterations can often be analyzed
		directly by verifying each $\ZF$ axiom within the particular
		construction, as in~\cite{Karagila2020Bristol}.
		The present framework supplies the limit-filter technology
		($\omega_1$-completeness, Theorem~\ref{thm:filter-generation}) that
		such a verification may use for ground-coded countable tuples of
		hereditarily symmetric names.
		We leave the formulation of a general class-length theorem to future
		work.
	\end{remark}
	
	\section{Withdrawn application}
	\label{sec:withdrawn}

	The previously advertised application to a family
	$\langle P_\alpha:\alpha<\kappa\rangle$ of two-element sets of reals
	without a choice function is withdrawn.  There are three independent
	obstructions.

	\paragraph{A ZF obstruction.}
	The space $2^\omega$ is linearly ordered lexicographically in $\ZF$.
	Every nonempty finite subset of a linearly ordered set has a unique least
	element.  Thus, for any set-indexed family $\langle P_i:i\in I\rangle$ of
	nonempty finite subsets of $2^\omega$, Replacement applied to
	$i\mapsto\min_{\mathrm{lex}}P_i$ produces a choice function.  In
	particular, the claimed family of pairs of reals with no choice function is
	impossible in $\ZF$.

	\paragraph{Hereditary symmetry.}
	Under the stated step filter generated by the whole swap group, the names for
	the two individual generic reals are exchanged by the nontrivial automorphism
	and are not symmetric.  Although the unordered pair name is fixed as a set,
	its constituents are not HS; therefore that pair name is not hereditarily
	symmetric.  The earlier argument checked symmetry of the outer name but not
	the hereditary clause.

	\paragraph{The same-generic argument.}
	Name invariance under an automorphism relates valuations in $G$ and in the
	transported generic $gG$.  It does not, without a fixed compatible condition
	or an equivalent forcing argument, permit one to apply $g$ internally to a
	selected value while keeping the same generic $G$.  The earlier
	no-choice-function proof made exactly this unsupported step.

	For the same reasons, the proposed finite-support proof of failure of $\DC$
	does not establish its conclusion: it uses the same non-HS constituents and
	the same invalid same-generic inference.  The finite-support model may have
	other choice-theoretic properties, but none is proved here.  All results and
	corollaries depending on the withdrawn general $\DC$ theorem or on this
	application are withdrawn.

	\section{Conclusion}

	The retained result is algebraic and deliberately limited.  Under the stated
	coherent direct-limit hypotheses, the head-pullback construction gives a
	normal $\omega_1$-complete limit filter: at uncountable cofinality all
	countably many generator stages and conjugating automorphisms are jointly
	bounded, while at countable cofinality the required completion is imposed by
	definition.  This entails closure of HS names under canonical tuples formed
	from ground-model countable sequences of HS names.

	The standard Symmetric Extension Theorem separately yields $\ZF$ whenever
	the stage data constitute a genuine symmetric system.  Neither that theorem
	nor filter completeness alone yields closure of the symmetric model under
	arbitrary countable sequences from $V[G]$.  Accordingly, no general
	dependent-choice preservation theorem is claimed.

	The unordered-pairs application is withdrawn both because its proposed names
	fail hereditary symmetry under the stated filter and because the advertised
	choice failure for pairs of reals is impossible in $\ZF$: lexicographic
	minima provide a canonical selector.  The corrected scope of the paper is
	therefore the limit-filter calculation and the qualified syntactic closure
	statement recorded above.

	\par\medskip
	\begin{acknowledgments}
		This work builds on the modern development of symmetric extensions and their iteration theory.
		We especially acknowledge Karagila's foundational treatment of finite-support symmetric iterations~\cite{karagila2019},
		and the standard forcing background in Jech~\cite{Jech03} and Kunen~\cite{Kunen2011}.
		We also thank the authors of the other works cited herein.
		Any remaining errors are the author's alone.
	\end{acknowledgments}
	
	\smallskip
	\noindent\textbf{AI tools.} Large language models (LLMs)
	were used as a writing aid for copy-editing, typesetting, and error-checking suggestions. All mathematical statements,
	proofs, and final editorial decisions were made by the author with whom final responsibility resides.

\end{document}